\documentclass[journal,transmag]{IEEEtran}

\usepackage{amsfonts}
\usepackage{amsmath}
\usepackage{amssymb}

\ifCLASSINFOpdf
   \usepackage[pdftex]{graphicx}
\else
\fi
\usepackage{amsmath}
\newtheorem{theorem}{Theorem}[section]

\begin{document}
%
\title{A Dynamic Analysis of Energy Storage with Renewable and\\ Diesel Generation using Volterra Equations}

\author{\IEEEauthorblockN{Denis Sidorov\IEEEauthorrefmark{1,3},~\IEEEmembership{Senior Member,~IEEE},
Ildar Muftahov\IEEEauthorrefmark{1,2},
Nikita Tomin\IEEEauthorrefmark{1},~\IEEEmembership{Member,~IEEE},
Dmitriy Karamov\IEEEauthorrefmark{1}, \\
Daniil Panasetsky\IEEEauthorrefmark{1},~\IEEEmembership{Member,~IEEE},
Aliona Dreglea\IEEEauthorrefmark{1,3},
Fang Liu\IEEEauthorrefmark{4},
and
Aoife Foley\IEEEauthorrefmark{5},~\IEEEmembership{Member,~IEEE}}
\IEEEauthorblockA{\IEEEauthorrefmark{1}Energy Systems Institute of Siberian Branch of Russian Academy of Sciences, Irkutsk 664033, Russia}
\IEEEauthorblockA{\IEEEauthorrefmark{2}Main Computing Center, JSC Russian Railways,  Irkutsk 664005, Russia}
\IEEEauthorblockA{\IEEEauthorrefmark{3}Baikal School of BRICS, Irkutsk National Research Technical University, Irkutsk 664074, Russia}
\IEEEauthorblockA{\IEEEauthorrefmark{4}School of Automation, Central South University, Changsha 410083, PRC}
\IEEEauthorblockA{\IEEEauthorrefmark{5}
School of Mechanical and Aerospace Engineering, Queen's University Belfast,
Belfast BT9 5AH, UK}
\thanks{Manuscript received .. 2019; revised... 2019.
This work was supported by the NSFC-RFBR Exchange Program under NSFC Grant No.~61911530132 and RFBR Grant No.~19-58-53011. It is also partly supported by RFBR Grant 18-31-00206. 
Corresponding author: Denis~Sidorov (e-mail:~dsidorov@isem.irk.ru).}}


\markboth{Journal of \LaTeX\ Class Files,~Vol.~14, No.~8, August~2019}%
{Shell \MakeLowercase{\textit{et al.}}: Bare Demo of IEEEtran.cls for IEEE Transactions on Magnetics Journals}
%



\IEEEtitleabstractindextext{%
\begin{abstract}
Energy storage systems will play a key role in the power system of the twenty first century considering the large penetrations of variable renewable energy, growth in transport electrification and decentralisation of heating loads. Therefore reliable real time methods to optimise energy storage, demand response and generation are vital for power system operations. This paper presents a concise review of battery energy storage and an example of battery modelling for renewable energy applications and second details an adaptive approach to solve this load levelling problem with storage. A dynamic evolutionary model based on the first kind  Volterra integral equation is used in both cases. A direct regularised numerical method is employed to find the least-cost dispatch of the battery in terms of integral equation  solution. Validation on real data shows that the proposed evolutionary Volterra model effectively generalises conventional discrete integral model taking into account both state of health and the availability of generation/storage.
\end{abstract}

}

\maketitle

\IEEEdisplaynontitleabstractindextext

%
\IEEEpeerreviewmaketitle

\section{Introduction}

\IEEEPARstart{T}{}he current development of energy systems requires  deep integration of renewable energy sources (RES) and energy storage systems (ESS) in both centralised and isolated  energy systems \cite{IRENA1}.  Such integration was observed for different power systems, starting with large backbone solar and wind farms, with small distributed generating complexes, and small installations located directly in  customer’s place \cite{IRENA2}. 
The stochasticity of wind and solar generation requires storage devices acting as  energy system stabilisers \cite{Cristobal3, Tsuanyo2015Modelingn5,  Merei2013Optimizationn4}. For large solar and wind farms, the pumped stations are usually used \cite{Christian2014Optimaln6,Benedikt2016Scenario-basedn7}. Such devices are related to medium and long-term storage and have an installed capacity from 10 MWh to 1000 MWh \cite{Sauer8}. 
In practice, electrochemical batteries are widely used. They allow instantly store and generate power to cover the load efficiently.  The  battery energy system connected to a centralised electrical network by means of a hybrid inverters is analysed in \cite{Dufo-Lopez2015Optimisationn9}. Such systems allow  to accumulate energy when the tariff is low and use it during the peak hours. Thus, it is possible to reduce the annual costs for the purchase of electrical energy. In \cite{Dufo-Lopez2016Optimisationn10}, the problem of optimising the installed capacity of batteries and solar panels in isolated power systems is attacked. In \cite{Dufo-Lopez2013Comparisonn11}, a comparison of various mathematical models predicting the service life of batteries in isolated energy systems is fulfilled. 
\cite{Ghosh} considered the photovoltaic--wind--diesel hybrid systems with hydrogen--based long--term storage, battery and diesel--generator. 
The hydrogen storage system is compared to the diesel-generator backup system to find the cost--effective system. 
It is found that the hydrogen-based storage is especially efficient for the higher latitudes depending on the seasonal renewable energy variation and fuel cost at the site of application.
In \cite{Dursun2011Comparativen12}, the problem of optimal control of batteries in an isolated energy system with renewable energy was considered. \cite{Stevens96} presents the results of a process for determining battery charging efficiency.
It is to be noted that the distributed energy storage (DES) units able to enhance the voltage stability and robustness of DC distribution networks using flexible voltage control \cite{Li19}. 
\cite{Li18}  proposed the coordinated control strategy to ensure a good performance of the frequency support under the wind variation and the  state of charge (SoC) changing.

The main tool for analysing  storage efficiency  is mathematical modeling. Depending on the problem specifics, the battery models must 
fulfill different, sometimes contradictory requirements. The linear models are often used to optimise the installed power of renewable energy sources and the capacity of storage devices. In \cite{L2005Designn17} the linear model of battery with constant efficiency is used. The coordinated control loops are designed for AC/DC shipboard power systems with the distributed energy storage in \cite{He2018An18}. The linear model of a battery  presented in \cite{Juan2012Optimumn19} takes into account the battery inverter efficiency
 depending on the load factor. The linear model  presented in \cite{Pablo2013Optimaln20} takes into account the service life and the cost of its disposal. All these models are based on determining the SoC of the battery with respect to each time interval $t \in [0, T].$ This takes into account technical limitations for charge and discharge. For example, lead-acid (OPzV, OPzS) batteries'  limitations are 20\% -25\% of the installed capacity. 
The Oliver Tremblay’s models implemented in Matlab / Simulink \cite{tremblay2009experimental22} allow to simulate lead-acid, lithium-ion, nickel-cadmium and nickel-metal hydride batteries. In addition, Tremblay’s models allow to take into account the nonlinear dependence of the no-load voltage on the SoC. The model  based on the generalised Shepherd ratio first presented in \cite{Shepherd1965EV23}. The Shepherd model contains a nonlinear term characterising the magnitude of the polarisation voltage depending on the current amplitude and the actual SoC of the battery.

In a real battery in idle mode, its voltage increases almost to the electromotive force (EMF) of idle, and when  discharge current appears, the voltage drops sharply. The presence of a nonlinear term allows to determine the actual discharge current of the battery, however, with a numerical solution, this leads to the algebraic cycle and makes the model unstable. In the Tremblay model, the voltage value is unambiguously determined by the values of the discharge current and the actual  battery charge level, which ensures sufficiently accurate simulation results for the discharge and charge modes of various types of batteries, including those used in energy systems with renewable energy sources.  The problem of optimal power control in isolated energy systems with renewable energy sources and batteries using the Tremblay model is considered  in \cite{Obukhov2017Sim25} and in other papers.
The active mass degradation of rechargeable batteries is analysed  in  \cite{Dirk2002Optimumn28}, \cite{Schiffer26} and in other publications. These and other models have been developed in the methodology for categorising batteries depending on operating conditions \cite{Svoboda2007Operatingn30} as part of the EU project  ``Benchmarking''. As result of this project, six key indices were proposed. The numerical analysis of these indices allows to accurately describe the processes of degradation of the active mass of batteries under the current operating conditions. Based on these indices it is possible to determine the battery type and hardware that is the most suitable for the specific operating conditions causing the degradation processes minimisation. This technology is widely used in the isolated energy systems with solar power plants.

The analysis  shows that the tasks related to modeling, optimisation of installed capacity, development of new structural materials and storage control strategies are among the key tasks at the present time. In most cases in system energy studies, batteries are represented by linear discrete models, the developments of new fundamental approaches describing energy storage processes is an important direction. Such studies should be fulfilled combining the solid fundamental results and practical experience. 

Fredholm and Volterra (evolutionary) integral equations are in the core of many mathematical models in
physics, energetics, economics and ecology. Volterra equations are among the classic approaches  in electro-analytical chemistry~\cite{bieniasz}.
Historical overview of the results concerning the Volterra integral equations
(VIEs) of the first kind is given in \cite{brunner}.
The theory of integral models of evolving systems was initiated  in the  works of  L. Kantorovich, R. Solow and V. Glushkov in the mid-20th Century.  
Such theory  employs  the VIEs of the first kind where bounds of the integration interval can be  functions of time. Here readers may refer to the monographs  \cite{hriton,apart, sidbook}.
These integral dynamical models take into account
the memory of a dynamical system when its past impacts its future evolution and can be employed  for  dynamic analysis of energy storage.

In present research the studies \cite{D2018Volterran31} are developed further in order to implement energy storage dynamic analysis using the evolutionary (Volterra) dynamical models. Dynamic system analysis is carried out on 
the conventional isolated electric power system  consisting of the photovoltaic arrays and  battery energy storage  operated in
parallel with
diesel generator backup system to serve the residential electric load. Here readers may refer to \cite{Ghosh}.

The paper is organised as follows. In Sec. II the new Volterra model of storage is presented, the connection of Volterra model with conventional ampere-hour integral model is established and formulated in terms of inverse and direct problems correspondingly. Sec. III presents the methodology of battery modeling. 
The results of the model's verification on the real dataset are discussed in Sec. IV.
The application of deep learning to dynamic analysis of microgrid with energy storage is considered in Sec.~V.
 Finally, Sec. VI  presents concluding remarks, further development and perspectives of proposed Volterra models of energy storage.   
  
\section{Volterra Model}

The Volterra equations describe the systems state evolution  and refer to the memory expressed by integral over a time interval in the past. In fact, ampere-hour integral model (direct problem)

\begin{equation}
SOC(t) = SOC(0) +\int_0^t \eta (\cdot) \,  i(\tau) d\tau
\label{eq1}
\end{equation}
can be formulated as inverse problem, i.e. the Volterra integral equation (VIE) of the first kind with respect to the instantaneous battery current $ i(\tau)$ (assumed positive for charge, negative for discharge). Here  $\eta(\cdot)$ is the battery Coulombic efficiency.  SOC can be expressed both in \% and ampere-hours (or kWh).  In general settings efficiency  is function of global time $t$, integration parameter $\tau,$ temperature and other parameters. Let us now recall the generic form of VIE  

\begin{equation}
\int_0^t K(t,\tau) \,  x(\tau) d\tau = f(t), \,\, t \in [0,T].
\label{eq2}
\end{equation}
Here kernel $K(t,\tau)$ and source function $f(t)$ are assumed to be known (up to some measurement errors), $x(\tau)$ is the desired alternating power function.  It is known that solutions to integral equations of the first kind can be unstable and this is a well known ill-posed inverse problem. 
The solution
of linear integral equations of the first kind is of course the classical problem and has been addressed
by numerous authors. But only few authors studied these equations with
discontinuous kernels $K(t,\tau)$ especially in nonlinear case. In general,
VIE of the first kind can be solved by reduction to equations of the second kind,  regularisation algorithms developed for Fredholm equations can be also applied as well as direct discretisation methods used in present paper.

The theory and regularised numerical methods to relieve the ill-posedness of the problem are employed in Sec. III as suggested in papers \cite{I2016Numericn32,D2013Onn33}. Proposed approach to dynamical 
analysis of energy storage is based on the solid mathematical theory \cite{sidbook} including the following existence and uniqness theorem.

\begin{theorem} \cite{I2016Numericn32}. \label{th1} 
Let the following conditions take place: $K_i(t,\tau),$ $G_i(\tau,x(\tau))$
are continuous, $i=\overline{1,n},$ $\alpha_i(t)$ have continuous derivatives
with respect to $t,$ $K_n(t,t)\neq 0,$ $\alpha_i(0)=0,$ $f(0)=0,$ $0<\alpha_1(t)< \dots <
\alpha_{n-1}(t) < t,$ $t \in (0,T).$ Let the functions $G_i(\tau, x(\tau))$
satisfy the Lipschitz conditions on $x$ with the Lipschitz constants $q_i $ such as
$$q_n + \sum\limits_{i=1}^{n-1} \alpha_i'(0) |K_n(0,0)^{-1}(K_i(0,0) - K_{i+1}(0,0))| (1+q_i) <1.$$
Then exists $\tau^*>0$ such that the Volterra equation
\begin{equation}
\int_0^t K(t,\tau, x(\tau)) \, d \tau = f(t)
\label{eq1_v}
\end{equation}
with jump discontinuous kernel 
  \begin{equation}
  K(t,\tau,x(\tau)) = \left\{ \begin{array}{ll}
         \mbox{$K_1(t,\tau)G_1(\tau,x(\tau)), \,\,\, t,\tau \in D_1$}, \\
         \mbox{\,\dots \,\,\,\,\, \dots \dots} \\
         \mbox{$K_n(t,\tau)G_n(\tau,x(\tau)), \,\,\, t,\tau \in D_n$}. \\
        \end{array} \right. 
\label{eq2_2}
\end{equation}
has an unique local solution in ${\mathcal C}_{[0,\tau^*]}.$ Here  
$D_i=\{ t,\tau | \alpha_{i-1}(t) < \tau < \alpha_i(t) \},$ $\alpha_0(t)=0,$ $\alpha_n(t)=t.$
Furthermore, if $\min\limits_{\tau^*\leq t \leq T} (t- \alpha_{n-1}(t)) = h>0$
then solution can be constructed on the whole interval $[\tau^*, T]$
using the step method combined with successive approximations.
Thereby equation (\ref{eq1_v}) has the unique global solution in 
${\mathcal C}_{[0,T]}.$ 
\end{theorem} 

Theorem~\ref{th1} guarantees the existence of  the unique alternating power function $x(t)$ based on the known storage systems' efficiency described in terms of the nonlinear kernel $K$  and the 
misbalance described by the source function $f(t)$. This methodology will be described next in details.

\section{Battery Modelling Methodology}

\begin{figure}[h!t]
\centering
\includegraphics[width=2.6in]{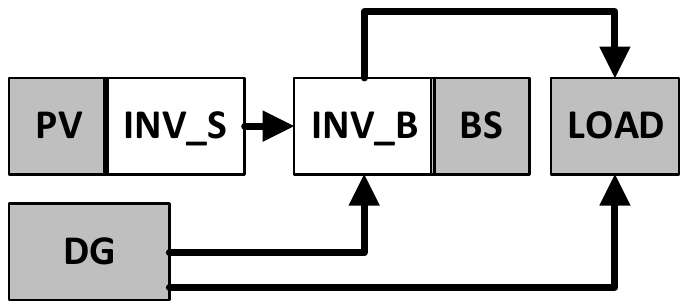}
\caption{The scheme of micrigrid with 
PV, battery energy storage system.
}
\label{fig1}
\end{figure}

The battery operation modes modeling is demonstrated using the example of an isolated hybrid energy system. It is assumed that this hybrid system uses: photovoltaic arrays (PV), solar inverters (INV\_S), battery inverters (INV\_B), battery energy storage  (BS),  and diesel generator (DG). 
Fig.~\ref{fig1} shows an isolated photovoltaic system  scheme used for the experiments and model validation below.

Installed capacities are as follows: PV is 75 kW, INV\_s is 75 kW, INV\_b is 72 kW, DG is 2x100 kW and BS is 384 kWh.
The maximum load is 47 kW.  Modeling the operation of a solar power station is performed using actinometric data recorded in the territory under consideration (see Sec. IV).

\subsection{Linear model}

Let us briefly describe the classic discrete  mathematical model with constant efficiency commonly used in the literature as discussed in Sec.~I. In fact, instead of the latter integral form (1), the following rather trivial linear  discrete model is used
\begin{equation}
SOC(t) = SOC(t-1)+I_s(t) \Delta t
\label{eq3}
\end{equation}
with constrain $I_s(t) \leq r_{BS} Q_{BS}^{max}$
where $SOC(t-1)$ (kW) is battery SOC in time $t-1,$
$I_s(t)$ is the alternating power function (kW);
$r_{BS}$ is technical restriction on the charge and 
discharge of battery (from 20\% till 40\%);
$Q_{BS}^{max}$ is installed battery capacity, kWh.
If the battery is charged then $I_s(t)$ is
multiplied by a constant efficiency of the battery 
and inverter. This value is assumed to be $0.8$ as suggested in \cite{Stevens96}.
$\Delta t $ is a discrete step (1 hour) to determine 
the amount of energy released to the battery. 
The term "linear model" is used for this classic model below. 

\subsection{Dynamical  model}
In order to efficiently model the storage operation, the following  integral dynamical model 
with constraints is employed
\begin{equation}\label{eq_volt}
\displaystyle\left\{ \begin{array}{ll}
\displaystyle\mbox{$\int_{0}^{t} K(t,\tau,x(\tau)) \, d\tau = f(t), \,\,\,   0 \leq \tau\leq t \leq T, $}\\
\displaystyle\mbox{${v(t) = \int_{0}^{t}{x(\tau) \,d\tau}, \; {\displaystyle\max_{t \in [0, T]}}{\;|v(t)|} \leq v_{max}}$},\\
\displaystyle\mbox{$E_{min}(t) \leq \int_{0}^{t}{v(\tau)\, d\tau}\leq E_{max}(t).$}
\end{array} \right.
\end{equation}
Here source function $ f(t) $ is the energy imbalance  defined as follows
$$f(t)=f_{PV}(t)  - f_{load}(t),$$
where $ f_ {PV}(t) $ is the PV generation,  and $ f_ {load} (t) $ is the electric load. 
This imbalance is supposed to be covered by battery storage 
 operated in
parallel with
diesel generator backup system.
 
In  integral model (\ref{eq_volt}), the alternating function of changing the power $ x(t) $ is the desired one. It allows for known maximum speed of the charge $ v_{max} $: 
\begin{enumerate}
\item
to determine $ E(t), $ which is the storage SoC under the constraints 
$ E_{min} (t) \leq E (t) \leq E_{max} (t) $ depending on the type of storage; 

\item to determine the minimum total capacity of the storage  to cover the load shortage; 

\item to calculate the number of  cycles based on behaviour of the  function $ E(t)$; 

\item the storage lifetime prediction  for the specific region.

\end{enumerate}

The problem of solution to the Volterra equation in problem (\ref{eq_volt}) is the typical inverse problem. Its integral operator is nonlinear in general  case and existence of the unique solution  follows from Theorem~\ref{th1}. Moreover, it is known (see Sec.~1.4 in \cite{apart})  that 
application of the discrete approximation methods to the
 Volterra equations of the first kind enjoys self-regularisation property when step-size is 
coordinated with amount of noise in the problem. The results of application of the regularisation procedure are discussed below in Sec. IV.
Following \cite{I2016Numericn32,D2013Onn33}, the error-resilient numerical method is employed.

The data flow diagram of developed software is shown in Fig.~6.
It consists of three main parts: the forecasting block, the block for VIE numerical solution to find the desired SoC and the diesel control block.
It is to be noted that for efficient application of the proposed Volterra model it is necessary to construct the accurate forecasts of the PV generation $ f_ {PV} $ and consumer electric load $ f_{load}(t)$. 
 For sake of simplicity PV generation  and  electric load are assumed to be known in this paper. In reality these time series are supposed to be forecasted
and therefore they are known up to some error level. This will be discussed in paragraph C, Sec. IV.

\begin{figure}[!t]
\centering
\includegraphics[width=3.2in]{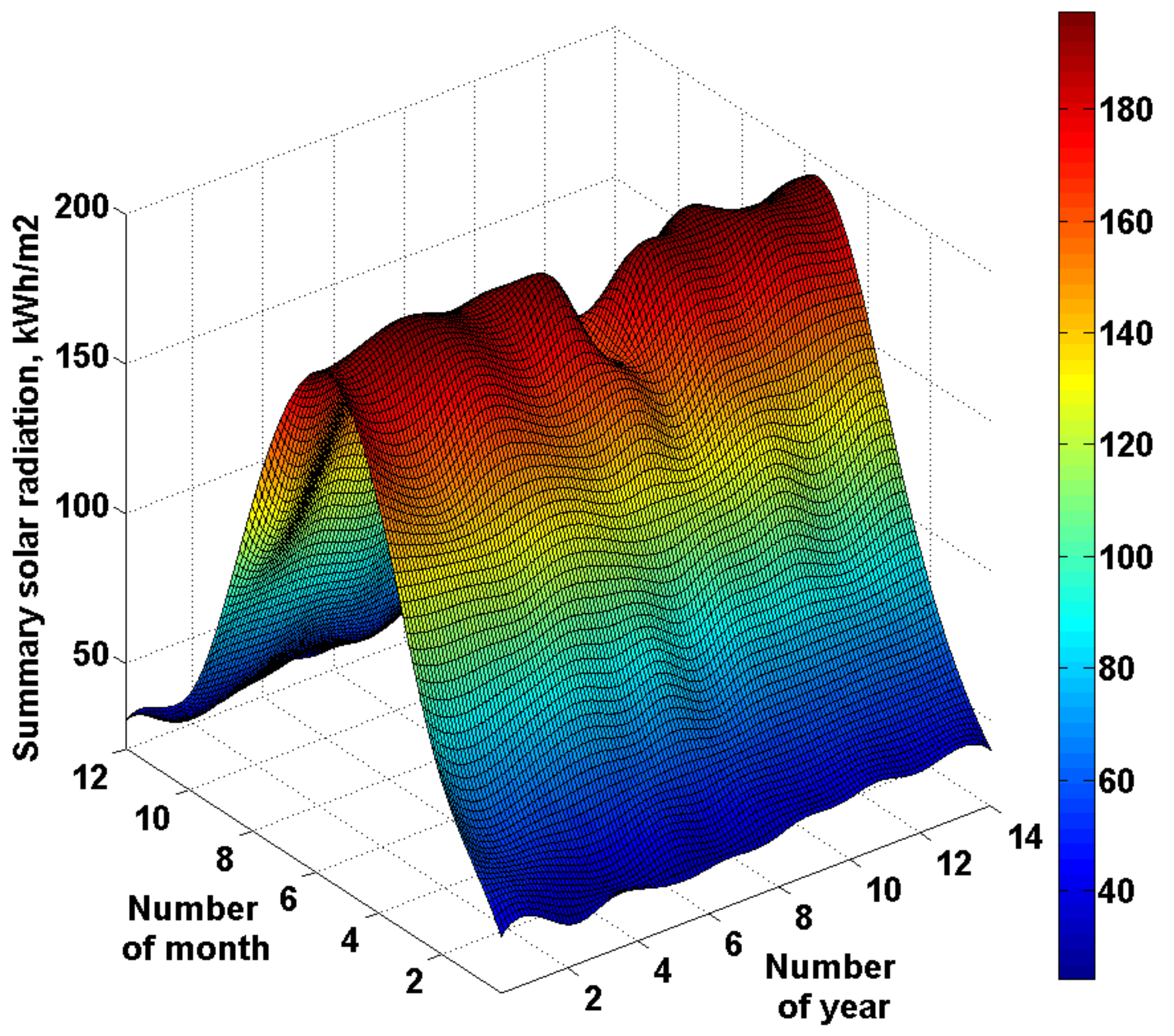}
\caption{Cumulative solar radiation for the considered region}
\label{fig2}
\end{figure}

\begin{figure}[!t]
\centering
\includegraphics[width=3.5in]{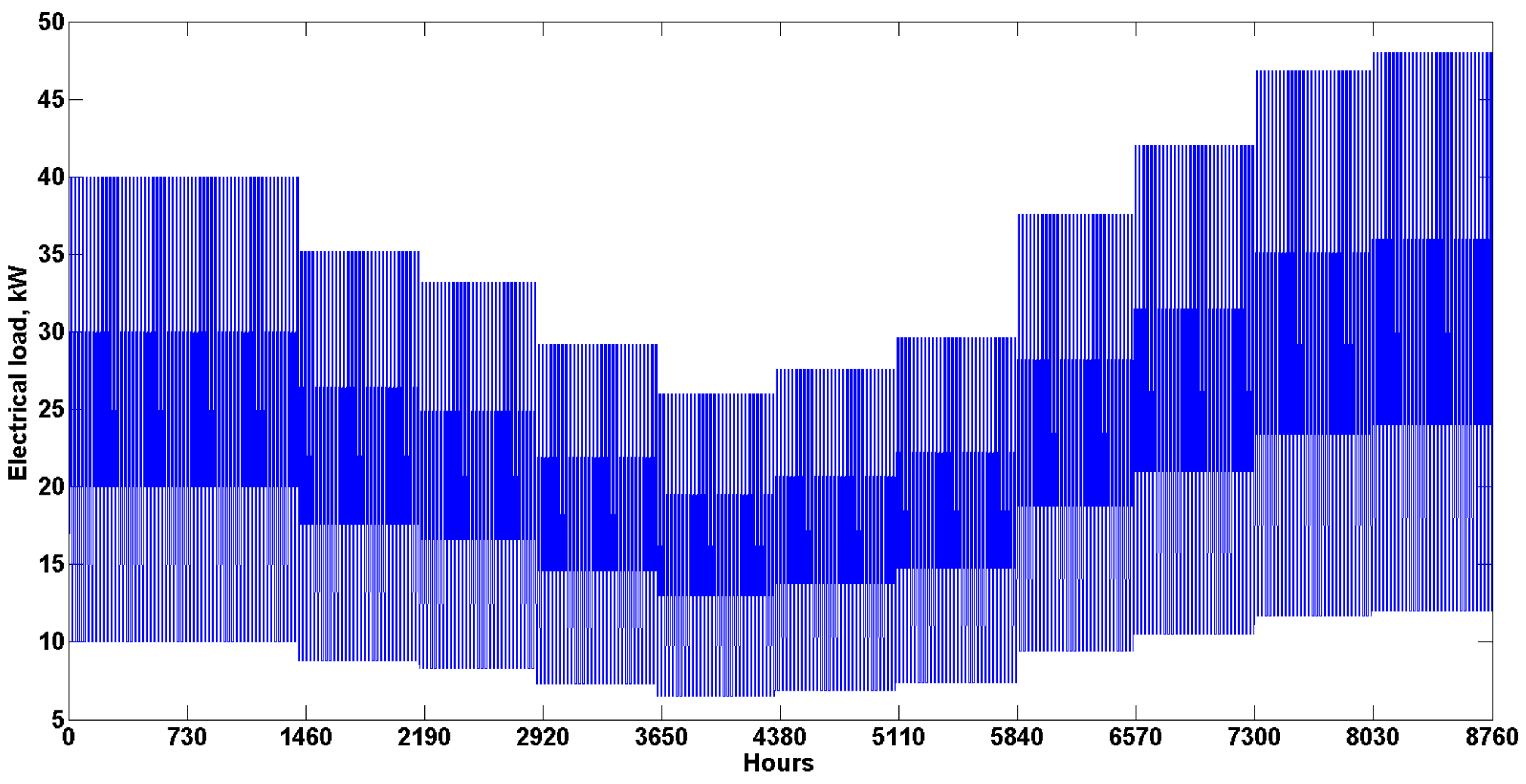}
\caption{Electrical load during the year}
\label{fig3}
\end{figure}

The load has maximum values in the winter and autumn periods as shown in Fig. \ref{fig3}.

 The following conditions are met:
\begin{enumerate}
\item if the alternating power function has a positive sign, then generation is sufficient for direct supply of the consumer and energy storage;
\item if the alternating power function has a negative sign, then energy is not enough. Therefore, the missing energy is taken from the battery;
\item if the battery SoC is below 20\%, then the diesel generator is turned on. The diesel generator turns on at full capacity covering the load and charging the battery.
\end{enumerate}

The system is simulated over the entire period of annual meteorological observations with a discrete step of one hour (a total of 122640 steps). After the simulation, the average annual value of the accumulated energy and monthly average values of SoC  are calculated.

Lead-carbon batteries adapted to heavy cyclic modes were used in the simulation. 
The study addresses two cases:
\begin{itemize}
\item case 1: the linear model of a battery with a constant efficiency.
\item case 2: the Volterra integral model with constant efficiency.
\end{itemize}

\paragraph{Case 1}
Simulation of an isolated energy system with the linear model of a battery with a constant efficiency \eqref{eq3} has the following average hourly values of  battery SoC.

\begin{figure}[!t]
\centering
\includegraphics[width=3.5in]{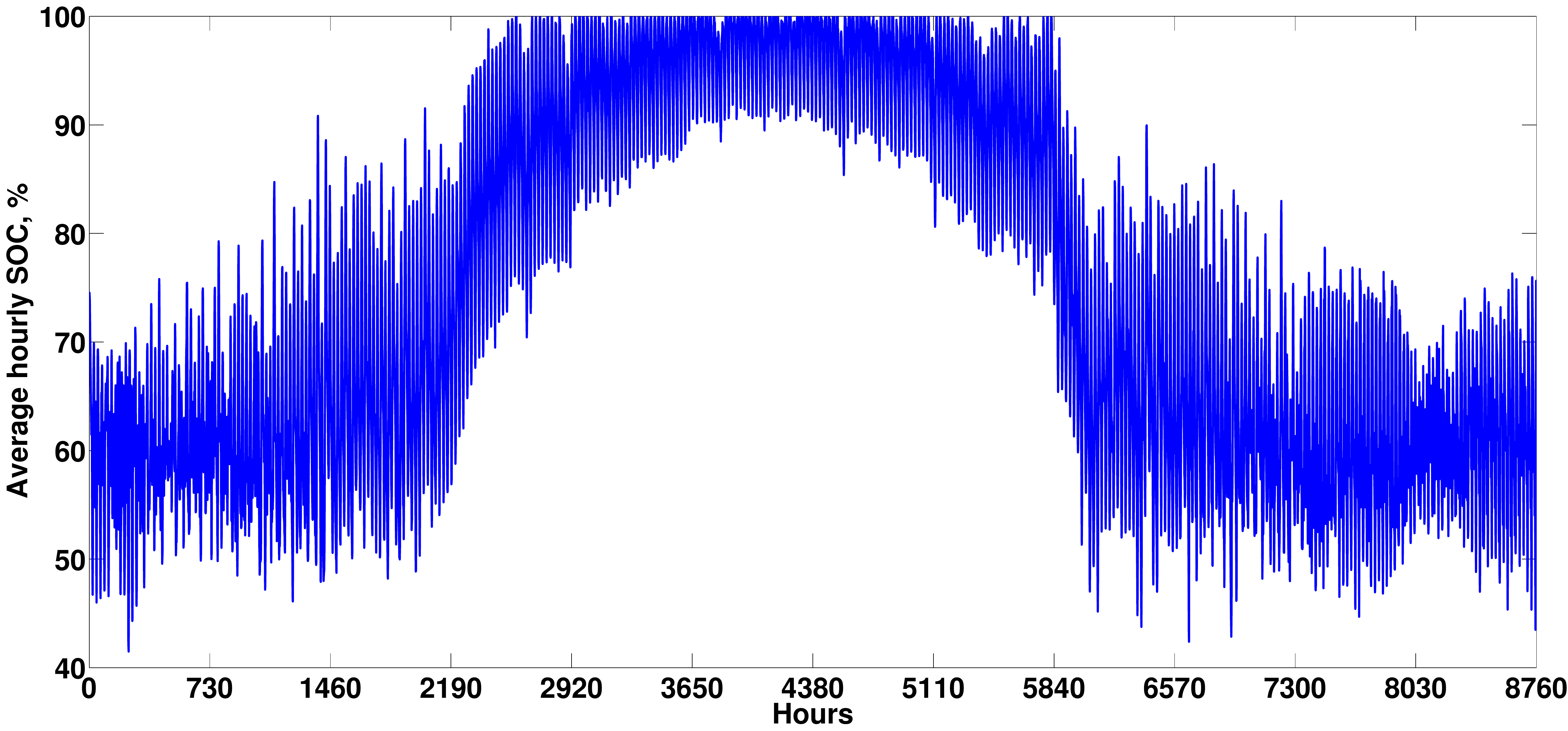}
\caption{Average SoC by the linear model}
\label{fig4}
\end{figure}

\begin{figure}[!t]
\centering
\includegraphics[width=3.5in]{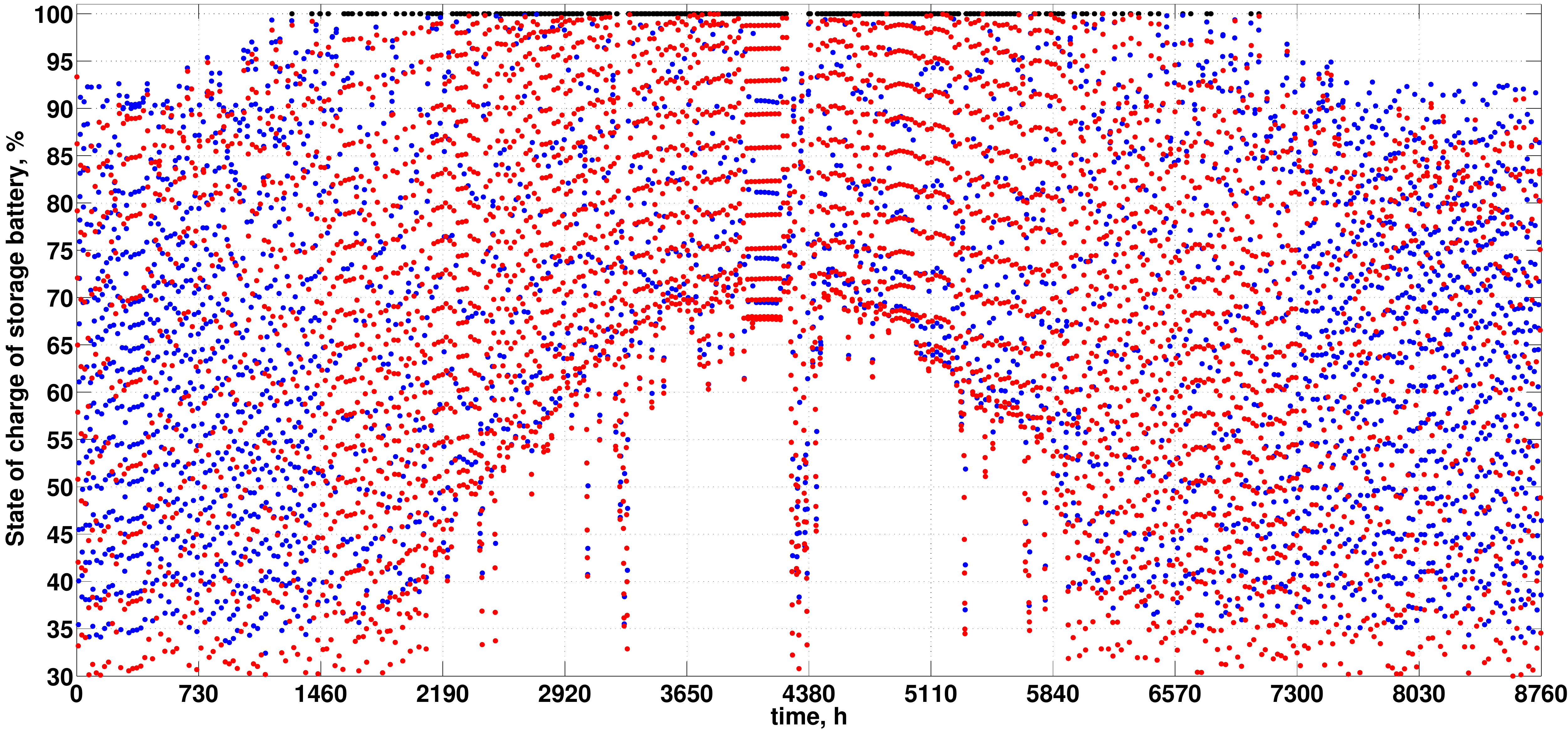}
\caption{SoC dynamics. Blue points shows charging, red stands for discharging and black is for idle state at high SoC }
\label{figSS}
\end{figure}

Calculations demonstrated that the average annual value of the energy supplied to the battery is 48116 kWh (PV is 20214 kWh, DG is 27902 kWh).
Fig.~\ref{fig4} shows the average SoC of the battery for the linear   model. Fig.~\ref{figSS}
demonstrates SoC dynamics over one-year period.

\paragraph{Case 2}
In the second case, the Volterra integral model with constant efficiency is used. 
By analogy with the first case,  the simulation of the system is performed throughout the entire period of meteorological observations, followed by averaging the necessary parameters.
Let us explain the meaning of the functions in model \eqref{eq2}. The right hand side consists of the difference between the PV generation and the load of consumers. Based on the unknown alternating power function $x(t)$ and restrictions on the maximum charge / discharge rate $v_{max}$, the charge level of the storage is determined, which also imposes threshold limits. The efficiency of  storage can depend on time and expressed by kernel $K(t,\tau)$. 

\begin{figure}[!t]
\centering
\includegraphics[width=3.5in]{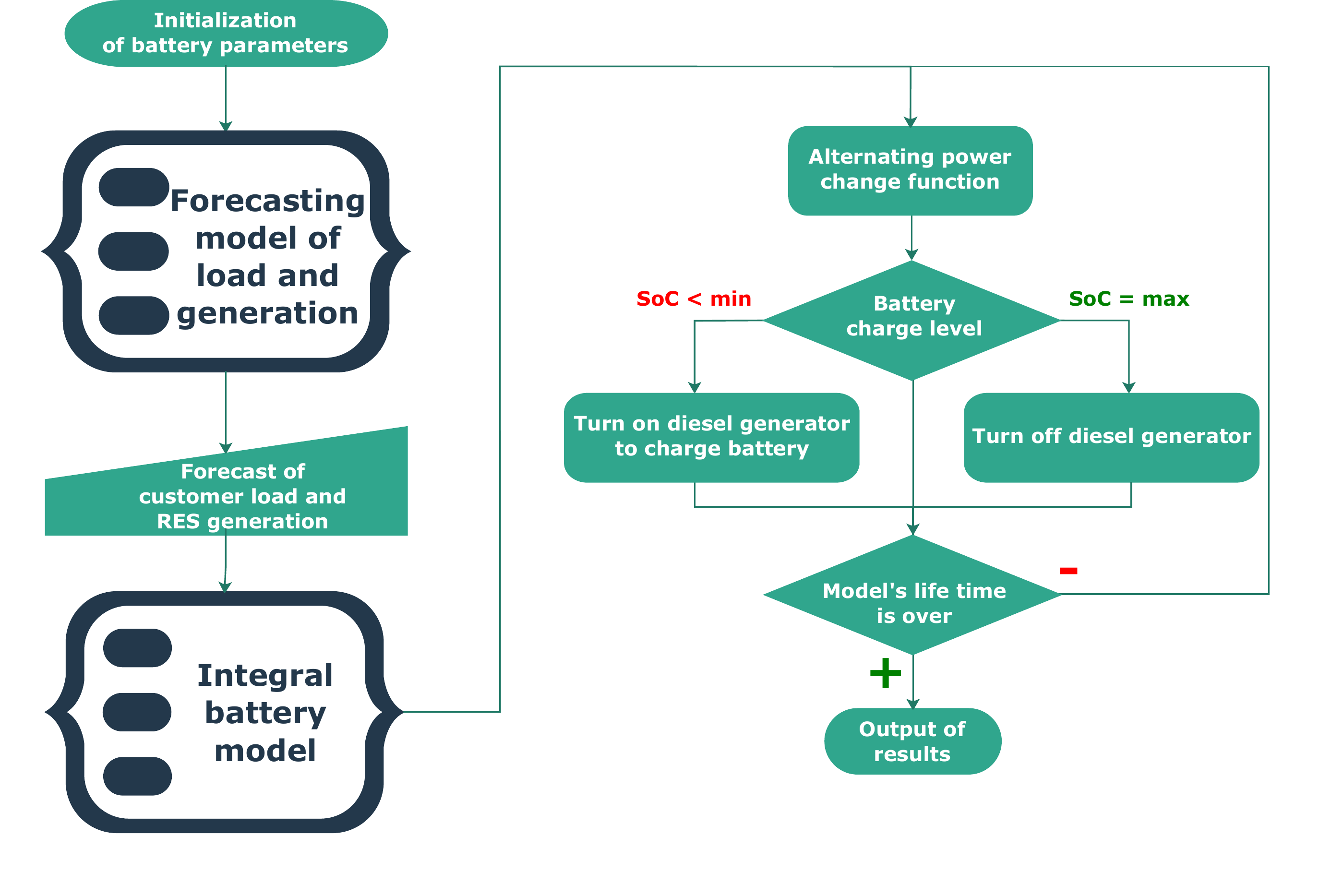}
\caption{Data flow diagram of the battery storage system }
\label{fig6}
\end{figure}

\begin{figure}[!t]
\centering
\includegraphics[width=3.5in]{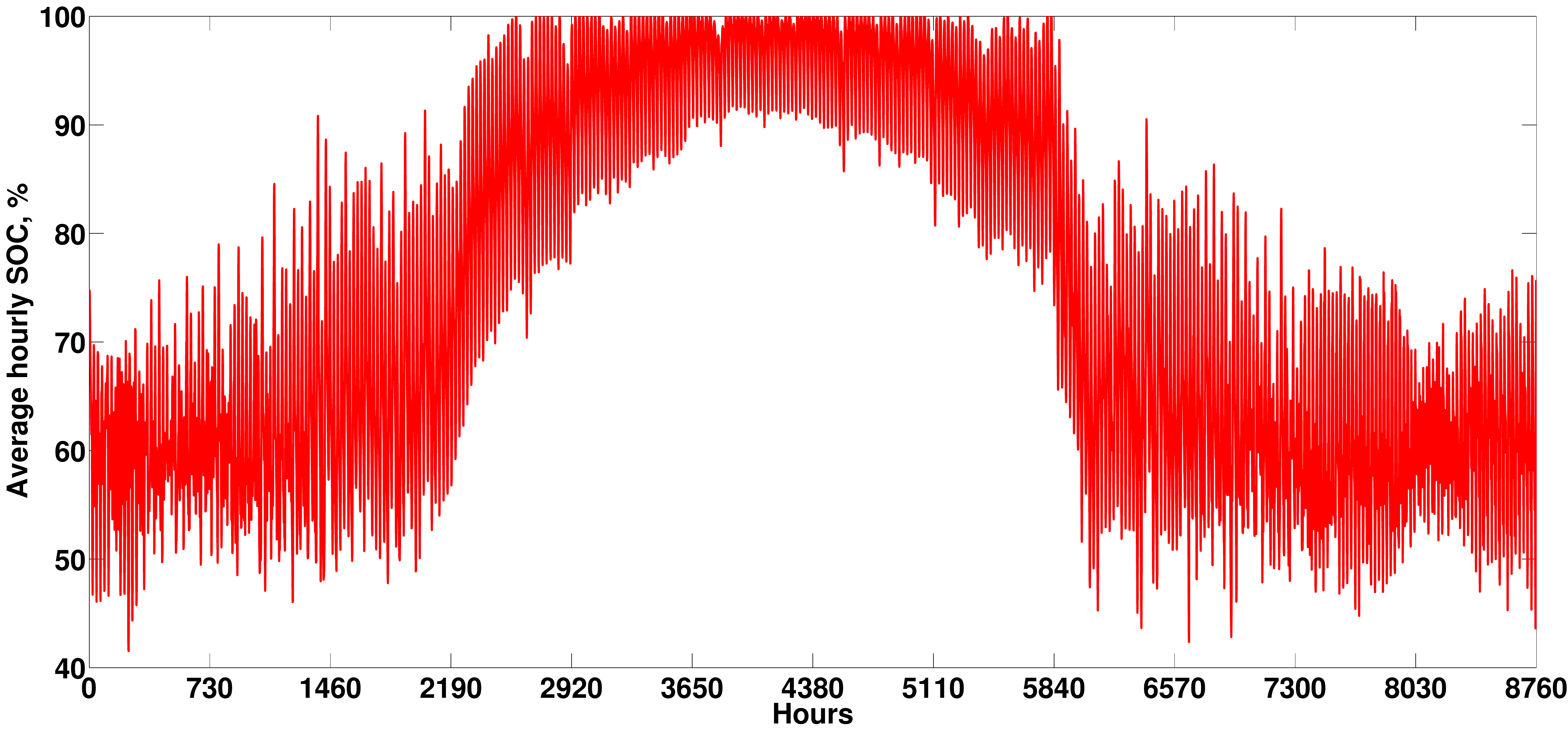}
\caption{Average SoC by the Volterra model}
\label{fig8}
\end{figure}

\begin{figure}[!t]
\centering
\includegraphics[width=3.5in]{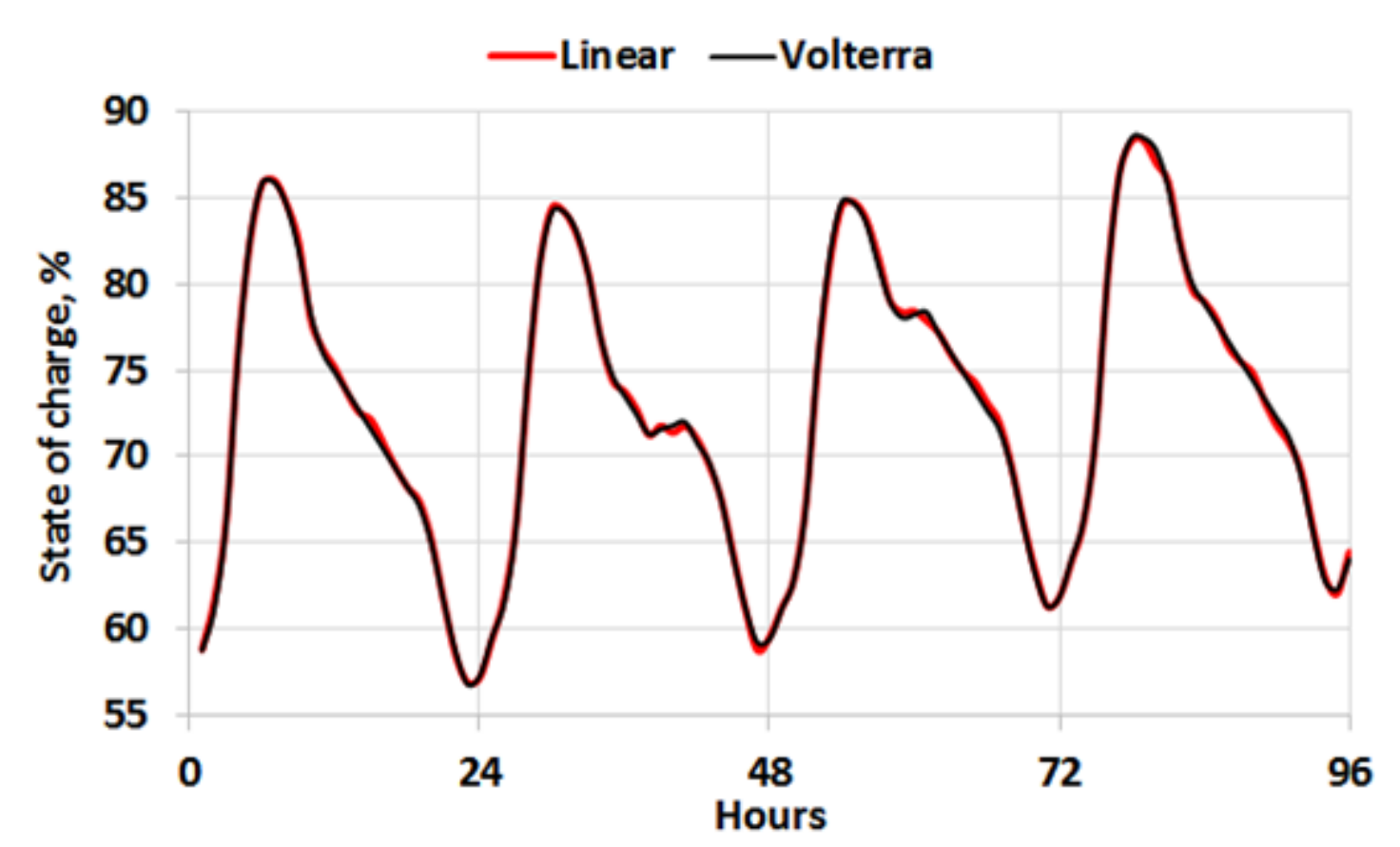}
\caption{SoC for conventional linear model and the Volterra model}
\label{fig9}
\end{figure}

Fig.~\ref{fig8} shows the average SoC of the battery for the Volterra model which is very similar with the results shown in Fig.~\ref{fig4}.
Both  SoC curves are shown in Fig.~\ref{fig9} for 96 hours period to demonstrate  the adequacy of the employed Volterra model. 
The results  show that the amount of energy supplied to the battery is 47849 kWh (PV: 20185 kWh, DG: 27664 kWh). Let us now consider the detailed validation of the Volterra model including the errors analysis and models robustness in input data.

\section{Model Validation}

\subsection{Test dataset description}
The isolated energy system of Innaly in Yakutia  was chosen for model validation. Retrospective data (time interval from 2005 to 2019) on the total solar radiation were used as baseline information. Solar radiation was recorded at a meteorological station located in this settlement. Electrical load  shown in Fig.~3 was built according to the real data of typical days relative to each month. Currently, this system has two diesel generators of 100 kW each. Their real technical characteristics were used in the work. The technical data of solar panels, PV and battery inverters, rechargeable batteries are used in the model validation below.  It is to be noted that storage efficiency is a non-linear function, it depends on the SoC, charging current and temperature. A detailed account of nonlinear dependencies is a difficult task to be described by Volterra kernel $K(t,\tau, x(\tau))$. In present study, the  efficiency was assumed to be constant ($\eta=0.8$) following \cite{Stevens96}.

\subsection{Model validation on the real dataset}

\begin{figure}[!t]
\centering
\includegraphics[width=3.5in]{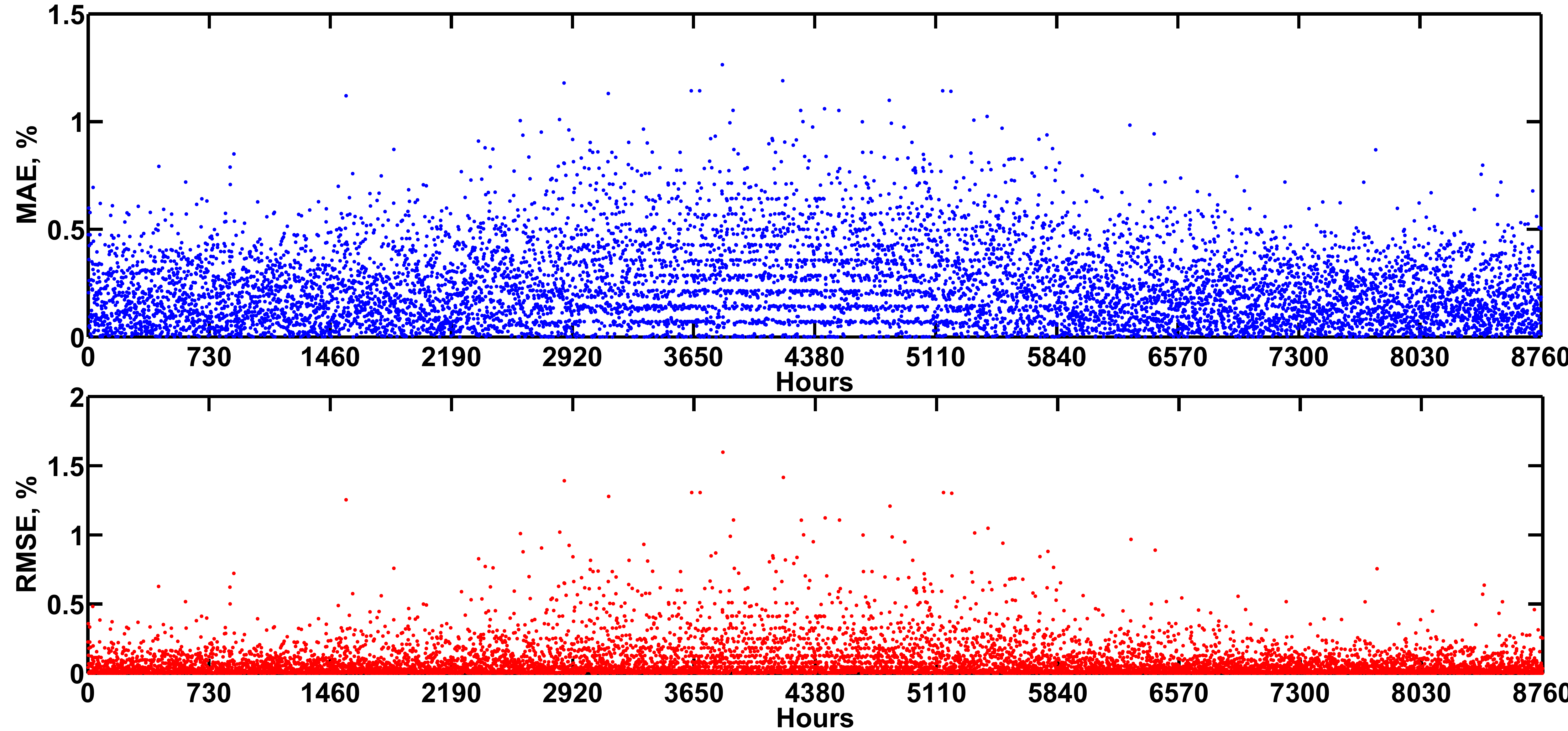}
\caption{RMSE and MAE for one year period}
\label{fig_new}
\end{figure}

The numerical experiments were fulfilled for the retrospective dataset.
  Three following metrics were selected:
$ RMSE ~ = ~ \sqrt {\frac{1}{n}\sum_{t = 1}^{n}(x_t-\bar{x_t})^2}$ is the root mean square error, mean absolute error $ MAE =\frac{1}{n}\sum_{t = 1}^{n}|x_t-\bar{x_t}|$ and the average absolute error in percent $MAPE=\frac{1}{n}\sum_ {t = 1}^{n}\frac {| x_t- \bar{x_t}|}{\bar{x_t}}*100\% $.
Here $\bar{x_t}$ corresponds to the linear model, and $x_t$ corresponds to the Volterra model. 
 The mean absolute error (MAE), root mean square error (RMSE) and mean absolute percentage error (MAPE) values are 0.23\%, 0.29\% and 0.31\%, respectively.
The results of multiple time scales are compared, the corresponding RMSE and MAE are shown in Fig.~\ref{fig_new}
 for the one year period. Each point corresponds to 1-hour sub-period.

 The simulation results of a battery using a linear model and Volterra integral model have demonstrated the similar results, see Fig. \ref{fig_new}. 

Thus, the model of SoC of the battery built on the Volterra equation enables the accurate processes description. The numerical results obtained in the calculation of the isolated energy system with PV, diesel and the battery showed the adequacy of Volterra model for the battery modeling. From the results it is clear that the results are very similar to the linear model. 
 
 \begin{figure}[!t]
\centering
\includegraphics[width=3.5in]{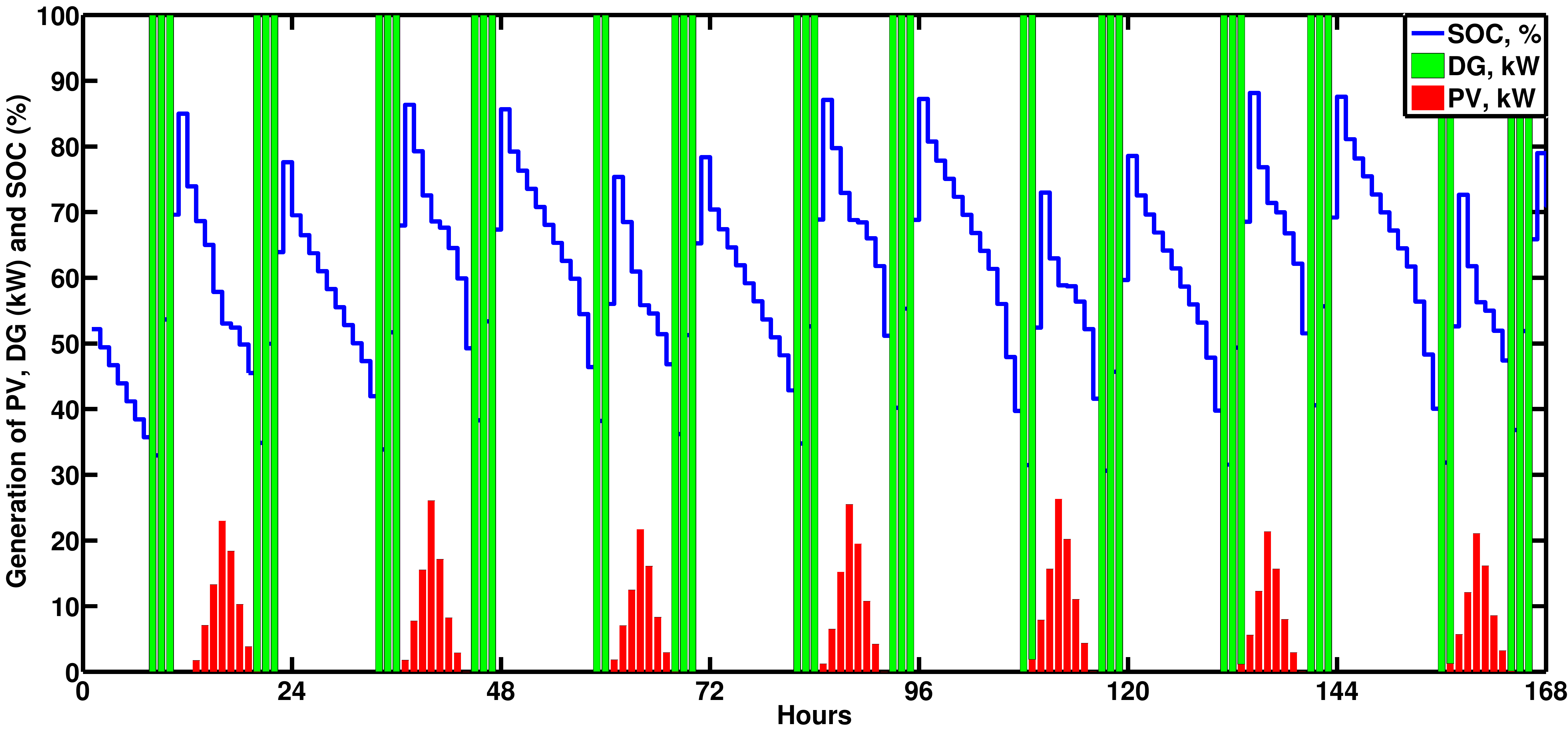}
\caption{Dynamic analysis of SoC of the battery  using the Volterra model with PV and diesel (January)}
\label{fig10}
\end{figure}

 \begin{figure}[!t]
\centering
\includegraphics[width=3.5in]{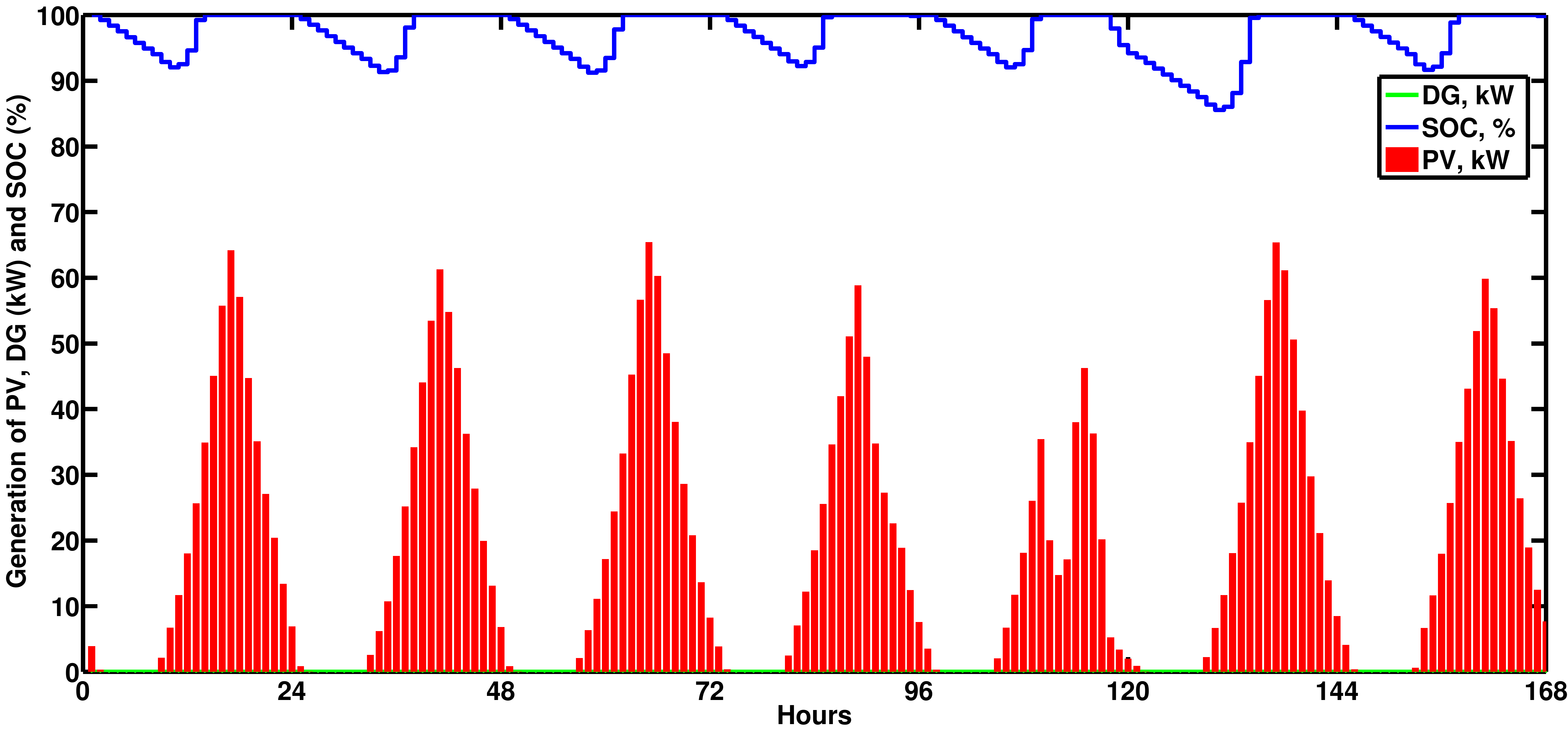}
\caption{Dynamic analysis of SoC of the battery  using the Volterra model with PV and diesel (June)}
\label{fig11}
\end{figure}

An analysis of Fig.~\ref{fig9} shows that the storage SoC repeats the behaviour of the alternating power function depending on the total solar radiation as shown in Fig. 2.
It should be noted that the Volterra integral equations have great potential in the form of accounting for nonlinear processes occurring in the battery. In this case, the consideration of nonlinearity is a characteristic feature inherent to Volterra integral equations. These processes include: nonlinear characteristics of efficiency, cycle life depending on the depth of discharge and processes of degradation of the active mass (SoH).

As footnote, Fig. \ref{fig10} and Fig. \ref{fig11} show typical dynamics of  storage SoC correspondingly for January and June. As can be noticed, proposed Volterra model with constrains can efficiently model the energy storage system with diesel  and PV.

Next paragraph focuses on the regularised Volterra models' robustness
to the noisy input data. 

\subsection{Regularised Volterra model}

One of the promising features of the Volterra model analysis of  storage is robustness
to the input data errors. As it is stated in the introduction, the problem of 
dynamic analysis of energy storage belongs to the class of ill-posed inverse problems.
It is to be noted that large amplifications of the measuring errors are typical for
ill-posed problems. 
In fact, the Volterra evolutionary dynamical models enjoys self-regularisation property when step-size is coordinated with amount of noise  in input data (or forecast error bounds). Here readers may refer to the monograph  \cite{apart}, p.~25
for more information concerning $(h, \alpha)$-regularisation.

But the step-size
is not under control in this problem
and therefore  the Lavrentiev $\alpha$-regularisation method is the only feasible option to attack this issue with solution stability.
Following  \cite{RegLavrISU2016}    the following regularised equation 
$$
	\alpha {x}_\alpha(t) + \int_0^t {K}(t, \tau) x_\alpha (\tau) d\tau = \tilde{f}(t)
$$ is used
 instead of the first equation in (\ref{eq_volt}). 
Here $\alpha$ is the regularisation parameter.
Following \cite{ccc2019} the publicly available dataset of the German power grid load, provided by ENTSO-E, for the period of 48 hours from 29.11.2013 22:00 till 01.12.2013 21:00 was used to demonstrate the superiority of regularised Volterra models when it comes to noisy input data. It is assumed that the exact load is unknown and only its forecast is available. In this experiment the storage efficiency (kernel ${K}(t, \tau)$) is assumed to  be constant, $\eta=0.92$  .

\begin{table}[t]

\renewcommand{\arraystretch}{1.3}
\caption{Errors analysis}
\label{tab1}
\centering

\begin{tabular}{lclclc|}
\hline

 &  GRU reg.   & GRU  \\
& $\alpha=0.422$ \\ 
\hline
 RMSE &  96.82    & 754.6  \\
 MAE &   75.12  & 537.88  \\
 MAPE&  16.78\%   &  130.05\%\\

\hline
 \end{tabular}
\end{table}

The deep learning GRU algorithm  (with settings: gru\_1(300) ,
 gru\_2(300), Dense\_1(16)) was used as base model for the electric load forecasting.  The  Lavrentiev $\alpha$-regularisation 
method is employed to cope with solution's  amplifications caused by inaccurate forecast.  
Table~\ref{tab1} demonstrates the errors analysis using RMSE, MAE and MAPE metrics when desired function is calculated based on the exact load, GRU-based load forecast without any regularisation and GRU-based load forecast with $\alpha$-regularisation. The regularisation parameter $\alpha=0.422$  is claimed here by the well-known  discrepancy principle. Regularisation significantly reduced the errors caused by inaccurate forecast.
It is to be noted that parameter $\alpha$ must be dynamically adjusted to forecast accuracy level.

As footnote, it can be concluded that $\alpha$-regularised Volterra
model is promising tool for mathematical modeling and dynamic analysis of power grids with energy storage.

\section{MDP-DQN Model}

  Optimal operation of a hybrid PV-diesel system can be formalised as a partially observable Markov decision process (MDP), where the hybrid system is considered as an agent that interacts with its environment \cite{Duan2019}. The fundamental difference between the approaches is that no specific strategy is set for the model, i.e. only the dynamics of a MDP environment and the conditions of the agent's actions. However, in the learning process, the agent finds the optimal policy (management strategy).

In order to approach the Markov property, the system's state $s_t\in S$ is made up of an history of features of observations $O^i_t$, $i \in {1, \ldots, N_f }$, where $N_f \in N$ is the total number of features. Each $O^i_t$ is represented by a sequence of punctual observations over a chosen history of length $h^i: O^i_t= [o^i_{t_{h^i+1}}, \ldots, o^i_t]$. At each time step, the agent observes a state variable $s_t$, takes an action $a_t \in A$ and moves into a state $s_{t+1}$. A reward signal $r_t = p(s_t, a_t, s_{t+1})$ is associated to the transition $(s_t, a_t, s_{t+1})$, where $p:\,S \times A \times S \rightarrow {\mathbb R} $ is the reward function. The  $\gamma$-discounted optimal Q-value function is defined as follows:
$$
Q^*(s,a) = \stackrel{max}{_{\pi}}\stackrel{E}{_{s_{t+1}, s_{t+2},\ldots}}\left[\sum^\infty_{k=t}
\gamma^{k-T} r_k | \,s_t=s, a_t=a, \pi\right].
$$
We propose to approximate $Q^*$ using a Deep Q-Network (DQN), because DQN models showed the good results for energy microgrids management in the recent works \cite{Mocanu2018, Xiao2018}. We adapted the approach, which was proposed in \cite{Francois2016} for a PV-hydrogen microgrid management. As result, we used DQN architecture, where the inputs are provided by the state vector, and where each separate output represents the Q-values for each discretiszed action. Possible actions are whether to turn on (or off) DG device for covering the load and charging the battery (avoid any value of loss load whenever possible). We considered three discretized actions : (i) turn on at full capacity, (ii) keep it idle or (iii) refill it, if DG tank is empty. 

\begin{figure}[h!t]
\centering
\includegraphics[width=3.5in]{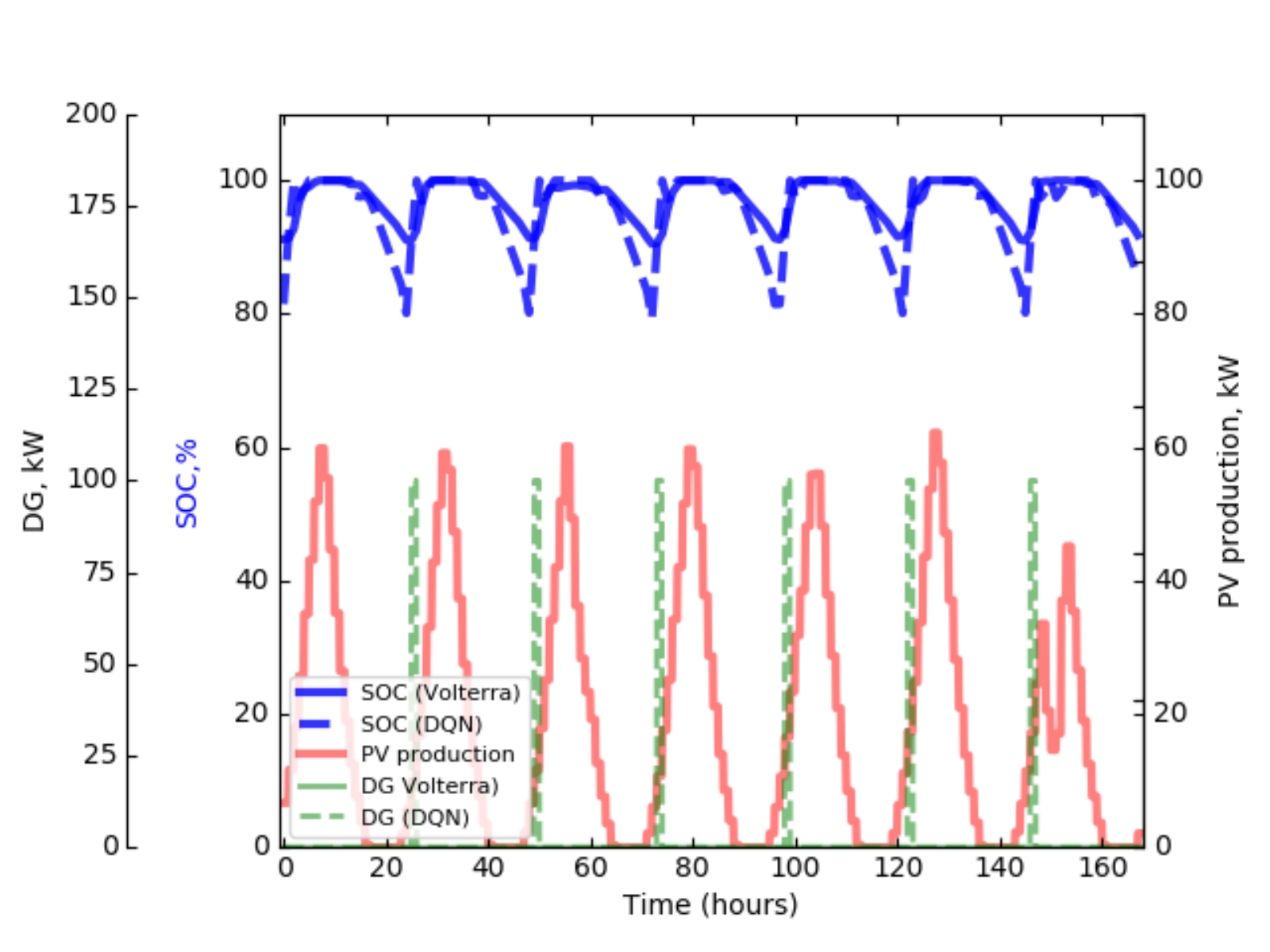}
\caption{Dynamic analysis of SoC of the battery using the Volterra and DQN models with PV and diesel (June)}
\label{figDQN}
\end{figure}

The reward function of the system corresponds to the instantaneous operational revenues $r_t$ at time $t \in T$. The instantaneous reward signal $r_t$ is obtained by adding the revenues
generated by the hydrogen production $r_DG$ with the penalties $r^-$ due to the value of loss load: 
$r_t = r(a_t, d_t) = r^{DG} (a_t, d_t) + r^-(a_t, d_t)$
, where $d_t$ denotes the net electricity demand. From the series of rewards ($r_t$), we obtain the operational revenues over year $y$ defined as follows: $M_y = \sum_{t\in\tau_y} r_t$  where $\tau_y$ is the set of time steps belonging to year $y$  \cite{Francois2016}.  
The typical behaviour of the policy for summer is illustrated in Figure \ref{figDQN}. This Figure illustrates the fact that the DQN model efficiently finds the policy, which is very similar to the policy (management strategy) of the Volterra model. 

\section{Conclusion \& Future Perspectives }

The objective of this paper is to fill the gap between the Volterra models and widely used storage models. An integral dynamical  model is proposed for  the  battery SoC dynamic analysis in terms of  the inverse problem's efficient solution. Proposed approach has the solid theoretical basis formulated in terms of the qualitative theory of integral equations and is complemented by the robust numerical methods development. The constructed linear Volterra model of the single battery storage with constant efficiency is verified on the real database and compared the results of the classic discrete model.

 The advantages of the suggested evolutionary dynamical model (\ref{eq1_v}) are as follows:
 \begin{enumerate}
 \item{regularised solution robust to  load and generation forecasts' unavoidable errors;}
 \item{definition of operating parameters of storage
when various renewable energy sources and storage are used jointly;}
\item{consideration of such characteristics of storage operation as power, charge/discharge rate, maximum number of work cycles, SoC limit;}
\item{low computational complexity of the algorithm when using a large number of storage units $n$;}
\item{accounting for the nonlinear nature of storages efficiency changes as function of SoC.}
 \end{enumerate}

The well developed theory of evolutionary linear and nonlinear integral equations (and their systems \cite{sidorov2013}) with continuous and jump discontinuous kernels opens new avenues for various storages models. Unlike the classic discrete models, proposed continuous models formulated in terms of integral equations suggests the error-resilient generic solutions. The Volterra evolutionary models naturally take into account the temporal degradation of the storages and their efficiency dependence of SoC.   The systems of Volterra equations can model energy storage systems combining the different  energy storage technologies.  The state of the art machine learning methods including DQN combined with the Volterra models will bring more insight to achieving the efficient balance of renewable sources and user demands in grids using the energy storages.

\section*{Acknowledgment}
We thank the anonymous reviewers for their helpful and constructive 
comments that greatly improved the manuscript.

Development of the mathematical theory  of this work was carried out under the State Assignment,  Project III.17.3 (reg. no. AAAA-A17-117030310442-8) and Project III.17.4 (reg. no. AAAA-A17-117030310438-1). 

This artilcle was finalised while first author was on a visiting research professor appointment at the School of Mechanical and Aerospace Engineering, Queen's University Belfast.

\ifCLASSOPTIONcaptionsoff
  \newpage
\fi



\bibliographystyle{IEEEtran}

\bibliography{bibfile}

\begin{thebibliography}{10}
\providecommand{\url}[1]{#1}
\csname url@samestyle\endcsname
\providecommand{\newblock}{\relax}
\providecommand{\bibinfo}[2]{#2}
\providecommand{\BIBentrySTDinterwordspacing}{\spaceskip=0pt\relax}
\providecommand{\BIBentryALTinterwordstretchfactor}{4}
\providecommand{\BIBentryALTinterwordspacing}{\spaceskip=\fontdimen2\font plus
\BIBentryALTinterwordstretchfactor\fontdimen3\font minus
  \fontdimen4\font\relax}
\providecommand{\BIBforeignlanguage}[2]{{%
\expandafter\ifx\csname l@#1\endcsname\relax
\typeout{** WARNING: IEEEtran.bst: No hyphenation pattern has been}%
\typeout{** loaded for the language `#1'. Using the pattern for}%
\typeout{** the default language instead.}%
\else
\language=\csname l@#1\endcsname
\fi
#2}}
\providecommand{\BIBdecl}{\relax}
\BIBdecl

\bibitem{IRENA1}
\BIBentryALTinterwordspacing
IRENA, ``Renewable energy statistics 2018,'' The International Renewable Energy
  Agency, {A}bu {D}habi, Tech. Rep., 2018. [Online]. Available:
  \url{https://www.irena.org/-/media/Files/IRENA/Agency/Publication/2018/Jul}
\BIBentrySTDinterwordspacing

\bibitem{IRENA2}
\BIBentryALTinterwordspacing
------, ``Off-grid renewable energy solutions: Global and regional status and
  trend,'' The International Renewable Energy Agency, {A}bu {D}habi, Tech.
  Rep., 2018. [Online]. Available:
  \url{https://www.irena.org/-/media/Files/IRENA/Agency/Publication/2018/Jul/IRENA\_Off-grid\_RE\_Solutions\_2018}
\BIBentrySTDinterwordspacing

\bibitem{Cristobal3}
\BIBentryALTinterwordspacing
I.~R. Cristóbal-Monreal and R.~Dufo-López, ``Optimisation of
  photovoltaic–diesel–battery stand-alone systems minimising system
  weight,'' \emph{Energy Conversion and Management}, vol. 119, pp. 279 -- 288,
  2016. [Online]. Available:
  \url{http://www.sciencedirect.com/science/article/pii/S0196890416303004}
\BIBentrySTDinterwordspacing

\bibitem{Tsuanyo2015Modelingn5}
\BIBentryALTinterwordspacing
D.~Tsuanyo, Y.~Azoumah, D.~Aussel, and P.~Neveu, ``Modeling and optimization of
  batteryless hybrid {PV} (photovoltaic)/diesel systems for off-grid
  applications,'' \emph{Energy}, vol.~86, pp. 152--163, jun 2015. [Online].
  Available: \url{https://doi.org/10.1016/j.energy.2015.03.128}
\BIBentrySTDinterwordspacing

\bibitem{Merei2013Optimizationn4}
\BIBentryALTinterwordspacing
G.~Merei, C.~Berger, and D.~U. Sauer, ``Optimization of an off-grid hybrid
  {PV}{\textendash}wind{\textendash}diesel system with different battery
  technologies using genetic algorithm,'' \emph{Solar Energy}, vol.~97, pp.
  460--473, nov 2013. [Online]. Available:
  \url{https://doi.org/10.1016/j.solener.2013.08.016}
\BIBentrySTDinterwordspacing

\bibitem{Christian2014Optimaln6}
\BIBentryALTinterwordspacing
C.~Bussar, M.~Moos, R.~Alvarez, P.~Wolf, T.~Thien, H.~Chen, Z.~Cai,
  M.~Leuthold, D.~U. Sauer, and A.~Moser, ``Optimal allocation and capacity of
  energy storage systems in a future european power system with 100{\%}
  renewable energy generation,'' \emph{Energy Procedia}, vol.~46, pp. 40--47,
  2014. [Online]. Available: \url{https://doi.org/10.1016/j.egypro.2014.01.156}
\BIBentrySTDinterwordspacing

\bibitem{Benedikt2016Scenario-basedn7}
\BIBentryALTinterwordspacing
B.~Lunz, P.~Stocker, S.~Eckstein, A.~Nebel, S.~Samadi, B.~Erlach,
  M.~Fischedick, P.~Elsner, and D.~U. Sauer, ``Scenario-based comparative
  assessment of potential future electricity systems {\textendash} a new
  methodological approach using germany in 2050 as an example,'' \emph{Applied
  Energy}, vol. 171, pp. 555--580, jun 2016. [Online]. Available:
  \url{https://doi.org/10.1016/j.apenergy.2016.03.087}
\BIBentrySTDinterwordspacing

\bibitem{Sauer8}
\BIBentryALTinterwordspacing
D.~U. Sauer, ``Chapter 2 - classification of storage systems,'' in
  \emph{Electrochemical Energy Storage for Renewable Sources and Grid
  Balancing}, P.~T. Moseley and J.~Garche, Eds.\hskip 1em plus 0.5em minus
  0.4em\relax Amsterdam: Elsevier, 2015, pp. 13 -- 21. [Online]. Available:
  \url{http://www.sciencedirect.com/science/article/pii/B9780444626165000024}
\BIBentrySTDinterwordspacing

\bibitem{Dufo-Lopez2015Optimisationn9}
\BIBentryALTinterwordspacing
R.~Dufo-L{\'{o}}pez, ``Optimisation of size and control of grid-connected
  storage under real time electricity pricing conditions,'' \emph{Applied
  Energy}, vol. 140, pp. 395--408, feb 2015. [Online]. Available:
  \url{https://doi.org/10.1016/j.apenergy.2014.12.012}
\BIBentrySTDinterwordspacing

\bibitem{Dufo-Lopez2016Optimisationn10}
\BIBentryALTinterwordspacing
R.~Dufo-L{\'{o}}pez, I.~R. Crist{\'{o}}bal-Monreal, and J.~M. Yusta,
  ``Optimisation of {PV}-wind-diesel-battery stand-alone systems to minimise
  cost and maximise human development index and job creation,'' \emph{Renewable
  Energy}, vol.~94, pp. 280--293, aug 2016. [Online]. Available:
  \url{https://doi.org/10.1016/j.renene.2016.03.065}
\BIBentrySTDinterwordspacing

\bibitem{Dufo-Lopez2013Comparisonn11}
\BIBentryALTinterwordspacing
R.~Dufo-L{\'{o}}pez, J.~M. Lujano-Rojas, and J.~L. Bernal-Agust{\'{\i}}n,
  ``Comparison of different lead{\textendash}acid battery lifetime prediction
  models for use in simulation of stand-alone photovoltaic systems,''
  \emph{Applied Energy}, vol. 115, pp. 242--253, feb 2014. [Online]. Available:
  \url{https://doi.org/10.1016/j.apenergy.2013.11.021}
\BIBentrySTDinterwordspacing

\bibitem{Ghosh}
\BIBentryALTinterwordspacing
P.~Ghosh, B.~Emonts, and D.~Stolten, ``Comparison of hydrogen storage with
  diesel-generator system in a pv–wec hybrid system,'' \emph{Solar Energy},
  vol.~75, no.~3, pp. 187 -- 198, 2003. [Online]. Available:
  \url{http://www.sciencedirect.com/science/article/pii/S0038092X03002792}
\BIBentrySTDinterwordspacing

\bibitem{Dursun2011Comparativen12}
\BIBentryALTinterwordspacing
E.~Dursun and O.~Kilic, ``Comparative evaluation of different power management
  strategies of a stand-alone {PV}/wind/{PEMFC} hybrid power system,''
  \emph{International Journal of Electrical Power {\&} Energy Systems},
  vol.~34, no.~1, pp. 81--89, jan 2012. [Online]. Available:
  \url{https://doi.org/10.1016/j.ijepes.2011.08.025}
\BIBentrySTDinterwordspacing

\bibitem{Stevens96}
J.~W. {Stevens} and G.~P. {Corey}, ``A study of lead-acid battery efficiency
  near top-of-charge and the impact on pv system design,'' in \emph{Conference
  Record of the Twenty Fifth IEEE Photovoltaic Specialists Conference - 1996},
  May 1996, pp. 1485--1488.

\bibitem{Li19}
Y.~{Li}, L.~{He}, F.~{Liu}, C.~{Li}, Y.~{Cao}, and M.~{Shahidehpour},
  ``Flexible voltage control strategy considering distributed energy storages
  for dc distribution network,'' \emph{IEEE Transactions on Smart Grid},
  vol.~10, no.~1, pp. 163--172, Jan 2019.

\bibitem{Li18}
\BIBentryALTinterwordspacing
Y.~Li, L.~He, F.~Liu, Y.~Tan, Y.~Cao, L.~Luo, and M.~Shahidehpour, ``A dynamic
  coordinated control strategy of wtg-es combined system for short-term
  frequency support,'' \emph{Renewable Energy}, vol. 119, pp. 1 -- 11, 2018.
  [Online]. Available:
  \url{http://www.sciencedirect.com/science/article/pii/S0960148117311655}
\BIBentrySTDinterwordspacing

\bibitem{L2005Designn17}
\BIBentryALTinterwordspacing
J.~L. Bernal-Agust{\'{\i}}n, R.~Dufo-L{\'{o}}pez, and D.~M. Rivas-Ascaso,
  ``Design of isolated hybrid systems minimizing costs and pollutant
  emissions,'' \emph{Renewable Energy}, vol.~31, no.~14, pp. 2227--2244, nov
  2006. [Online]. Available: \url{https://doi.org/10.1016/j.renene.2005.11.002}
\BIBentrySTDinterwordspacing

\bibitem{He2018An18}
\BIBentryALTinterwordspacing
L.~He, Y.~Li, Z.~Shuai, J.~M. Guerrero, Y.~Cao, M.~Wen, W.~Wang, and J.~Shi,
  ``A flexible power control strategy for hybrid {AC}/{DC} zones of shipboard
  power system with distributed energy storages,'' \emph{{IEEE} Transactions on
  Industrial Informatics}, vol.~14, no.~12, pp. 5496--5508, dec 2018. [Online].
  Available: \url{https://doi.org/10.1109/tii.2018.2849201}
\BIBentrySTDinterwordspacing

\bibitem{Juan2012Optimumn19}
\BIBentryALTinterwordspacing
J.~M. Lujano-Rojas, C.~Monteiro, R.~Dufo-L{\'{o}}pez, and J.~L.
  Bernal-Agust{\'{\i}}n, ``Optimum load management strategy for
  wind/diesel/battery hybrid power systems,'' \emph{Renewable Energy}, vol.~44,
  pp. 288--295, aug 2012. [Online]. Available:
  \url{https://doi.org/10.1016/j.renene.2012.01.097}
\BIBentrySTDinterwordspacing

\bibitem{Pablo2013Optimaln20}
\BIBentryALTinterwordspacing
P.~Garc{\'{\i}}a, J.~P. Torreglosa, L.~M. Fern{\'{a}}ndez, and F.~Jurado,
  ``Optimal energy management system for stand-alone wind
  turbine/photovoltaic/hydrogen/battery hybrid system with supervisory control
  based on fuzzy logic,'' \emph{International Journal of Hydrogen Energy},
  vol.~38, no.~33, pp. 14\,146--14\,158, nov 2013. [Online]. Available:
  \url{https://doi.org/10.1016/j.ijhydene.2013.08.106}
\BIBentrySTDinterwordspacing

\bibitem{tremblay2009experimental22}
O.~Tremblay and L.-A. Dessaint, ``Experimental validation of a battery dynamic
  model for ev applications,'' \emph{World Electric Vehicle Journal}, vol.~3,
  no.~2, pp. 289--298, 2009.

\bibitem{Shepherd1965EV23}
C.~Shepherd, ``Design of primary and secondary cells. p. 2. an equation
  describing battery discharge,'' \emph{Journal of Electrochemical Society},
  vol.~2, pp. 657--664, 1965.

\bibitem{Obukhov2017Sim25}
S.~Obukhov and I.~Plotnikov, ``Simulation model of operation of autonomuus
  photovoltaic plant under actual operation conditions,'' \emph{Geo Assets
  Engineerings}, vol. 328, no.~6, pp. 38--51, 2017.

\bibitem{Dirk2002Optimumn28}
D.~U. Sauer and J.~Garche, ``Optimum battery design for applications in
  photovoltaic systems {\textemdash} theoretical considerations,''
  \emph{Journal of Power Sources}, vol.~95, no. 1-2, pp. 130--134, mar 2001.

\bibitem{Schiffer26}
\BIBentryALTinterwordspacing
J.~Schiffer, D.~U. Sauer, H.~Bindner, T.~Cronin, P.~Lundsager, and R.~Kaiser,
  ``Model prediction for ranking lead-acid batteries according to expected
  lifetime in renewable energy systems and autonomous power-supply systems,''
  \emph{Journal of Power Sources}, vol. 168, no.~1, pp. 66--78, 2007, 10th
  European Lead Battery Conference. [Online]. Available:
  \url{http://www.sciencedirect.com/science/article/pii/S0378775306025122}
\BIBentrySTDinterwordspacing

\bibitem{Svoboda2007Operatingn30}
\BIBentryALTinterwordspacing
V.~Svoboda, H.~Wenzl, R.~Kaiser, A.~Jossen, I.~Baring-Gould, J.~Manwell,
  P.~Lundsager, H.~Bindner, T.~Cronin, P.~N{\o}rg{\aa}rd, A.~Ruddell,
  A.~Perujo, K.~Douglas, C.~Rodrigues, A.~Joyce, S.~Tselepis, N.~van~der Borg,
  F.~Nieuwenhout, N.~Wilmot, F.~Mattera, and D.~U. Sauer, ``Operating
  conditions of batteries in off-grid renewable energy systems,'' \emph{Solar
  Energy}, vol.~81, no.~11, pp. 1409--1425, nov 2007. [Online]. Available:
  \url{https://doi.org/10.1016/j.solener.2006.12.009}
\BIBentrySTDinterwordspacing

\bibitem{bieniasz}
L.~Bieniasz, \emph{Modelling Electroanalytical Experiments by the Integral
  Equation Method}.\hskip 1em plus 0.5em minus 0.4em\relax Berlin: Springer,
  2015.

\bibitem{brunner}
\BIBentryALTinterwordspacing
H.~Brunner, ``1896–1996: One hundred years of {V}olterra integral equations
  of the first kind,'' \emph{Applied Numerical Mathematics}, vol.~24, no.~2,
  pp. 83 -- 93, 1997, second International Conference on the Numerical Solution
  of {V}olterra and Delay Equations. [Online]. Available:
  \url{http://www.sciencedirect.com/science/article/pii/S0168927497000135}
\BIBentrySTDinterwordspacing

\bibitem{hriton}
N.~Hritonenko and Y.~Yatsenko, \emph{Modeling and Optimization of the Lifetime
  of Technologies}.\hskip 1em plus 0.5em minus 0.4em\relax Dordrecht: Kluwer
  Academic Publishers, 1996.

\bibitem{apart}
A.~S. Apartsyn, \emph{Nonclassical linear {V}olterra equations of the first
  kind}.\hskip 1em plus 0.5em minus 0.4em\relax Walter de Gruyter, 2011,
  vol.~39.

\bibitem{sidbook}
\BIBentryALTinterwordspacing
D.~Sidorov, \emph{Integral Dynamical Models: Singularities, Signals and
  Control}.\hskip 1em plus 0.5em minus 0.4em\relax World Scientific Publ.,
  2014. [Online]. Available:
  \url{https://www.worldscientific.com/doi/abs/10.1142/9278}
\BIBentrySTDinterwordspacing

\bibitem{D2018Volterran31}
\BIBentryALTinterwordspacing
D.~Sidorov, A.~Zhukov, A.~Foley, A.~Tynda, I.~Muftahov, D.~Panasetsky, and
  Y.~Li, ``{V}olterra models in load leveling problem,'' \emph{E3S Web of
  Conferences}, vol.~69, p. 01015, 2018. [Online]. Available:
  \url{https://doi.org/10.1051/e3sconf/20186901015}
\BIBentrySTDinterwordspacing

\bibitem{I2016Numericn32}
\BIBentryALTinterwordspacing
I.~Muftahov, A.~Tynda, and D.~Sidorov, ``Numeric solution of {V}olterra
  integral equations of the first kind with discontinuous kernels,''
  \emph{Journal of Computational and Applied Mathematics}, vol. 313, pp.
  119--128, mar 2017. [Online]. Available:
  \url{https://doi.org/10.1016/j.cam.2016.09.003}
\BIBentrySTDinterwordspacing

\bibitem{D2013Onn33}
\BIBentryALTinterwordspacing
D.~N. Sidorov, ``On parametric families of solutions of {V}olterra integral
  equations of the first kind with piecewise smooth kernel,''
  \emph{Differential Equations}, vol.~49, no.~2, pp. 210--216, feb 2013.
  [Online]. Available: \url{https://doi.org/10.1134/s0012266113020079}
\BIBentrySTDinterwordspacing

\bibitem{RegLavrISU2016}
I.~R. Muftahov, D.~N. Sidorov, and N.~A. Sidorov, ``Lavrentiev regularization
  of integral equations of the first kind in the space of continuous
  functions,'' \emph{The Bulletin of Irkutsk State University. Series
  Mathematics}, vol.~15, pp. 62--77, 2016.

\bibitem{ccc2019}
D.~{Sidorov}, Q.~{Tao}, I.~{Muftahov}, A.~{Zhukov}, D.~{Karamov}, A.~{Dreglea},
  and F.~{Liu}, ``Energy balancing using charge/discharge storages control and
  load forecasts in a renewable-energy-based grids,'' in \emph{2019 38th
  Chinese Control Conference (CCC)}, vol. (to appear), no. arXiv:1906.02959
  [eess.SP], July 2019, pp. 1--6.

\bibitem{Duan2019}
J.~{Duan}, Z.~{Yi}, D.~{Shi}, C.~{Lin}, X.~{Lu}, and Z.~{Wang},
  ``Reinforcement-learning-based optimal control for hybrid energy storage
  systems in hybrid ac/dc microgrids,'' \emph{IEEE Transactions on Industrial
  Informatics}, pp. 1--1, 2019.

\bibitem{Mocanu2018}
E.~{Mocanu}, D.~C. {Mocanu}, P.~H. {Nguyen}, A.~{Liotta}, M.~E. {Webber},
  M.~{Gibescu}, and J.~G. {Slootweg}, ``On-line building energy optimization
  using deep reinforcement learning,'' \emph{IEEE Transactions on Smart Grid},
  pp. 1--1, 2018.

\bibitem{Xiao2018}
\BIBentryALTinterwordspacing
L.~Xiao, X.~Xiao, C.~Dai, M.~Peng, L.~Wang, and H.~V. Poor, ``Reinforcement
  learning-based energy trading for microgrids,'' \emph{CoRR}, vol.
  abs/1801.06285, 2018. [Online]. Available:
  \url{http://arxiv.org/abs/1801.06285}
\BIBentrySTDinterwordspacing

\bibitem{Francois2016}
V.~Francois-Lavet, Q.~Gemine, D.~Ernst, and R.~Fonteneau, ``Towards the
  minimization of the levelized energy costs of microgrids using both long-term
  and short-term storage devices.'' \emph{Smart Grid: Networking, Data
  Management, and Business Models}, pp. 295--319, 2016.

\bibitem{sidorov2013}
\BIBentryALTinterwordspacing
D.~N. Sidorov, ``Solvability of systems of {V}olterra integral equations of the
  first kind with piecewise continuous kernels,'' \emph{Russian Mathematics},
  vol.~57, no.~1, pp. 54--63, Jan 2013. [Online]. Available:
  \url{https://doi.org/10.3103/S1066369X13010064}
\BIBentrySTDinterwordspacing

\end{thebibliography}

\end{document}